\documentclass[journal]{IEEEtran}
\usepackage{cite}
\usepackage{amsmath,amsfonts,amssymb,amsthm}
\theoremstyle{definition}

\theoremstyle{plain}

\newtheorem{lemma}{Lemma}
\newtheorem{proposition}{Proposition}

\usepackage{graphicx}
\usepackage{mathrsfs}
\usepackage{booktabs}
\usepackage{multirow}
\usepackage{array}
\usepackage[ruled]{algorithm}
\usepackage{hyperref}
\hypersetup{colorlinks,%
	citecolor={red}, 
	urlcolor={blue},
	linkcolor={blue},
	breaklinks={true},
	pagebackref={true},
	hyperindex={true},
}
\usepackage{url}
\usepackage{xcolor}
\usepackage{pgf,tikz}
\usetikzlibrary{arrows,shapes,positioning}
\usepackage{circuitikz}
\usepackage{subfig}
\usepackage{setspace}
\usepackage{url}
\usepackage{algorithmic}
\floatname{algorithm}{Model}

\usepackage{enumitem}
\newlist{abbrv}{itemize}{1}
\setlist[abbrv,1]{label=,labelwidth=0.9in,align=parleft,noitemsep,leftmargin=!}

\newcommand{\R}{\mathbb{R}}
\newcommand{\C}{\mathbb{C}}

\newcommand{\Z}{\mathbb{Z}}

\newcommand{\Herm}{\mathbb{H}}
\newcommand{\G}{\mathscr{G}}

\newcommand{\E}{\mathcal{E}}
\newcommand{\K}{\mathcal{K}}
\newcommand{\cplx}[1]{\mathrm{#1}}
\newcommand{\rv}[1]{\boldsymbol{#1}}
\newcommand{\cv}[1]{\boldsymbol{\mathrm{#1}}}
\newcommand{\ub}[1]{\overline{#1}}
\newcommand{\lb}[1]{\underline{#1}}
\newcommand{\Node}{\mathcal{N}}
\newcommand{\Gen}{\mathcal{G}}
\newcommand{\Shunt}{\mathcal{U}}
\newcommand{\Branch}{\mathcal{L}}
\newcommand{\Tap}{\mathcal{T}}
\newcommand{\Reseau}{\mathscr{P}}

\newcommand{\jc}{\cplx{j}}
\newcommand{\pl}{\lb{p}}
\newcommand{\ql}{\lb{q}}
\newcommand{\pu}{\ub{p}}
\newcommand{\qu}{\ub{q}}
\newcommand{\vl}{\lb{v}}
\newcommand{\vu}{\ub{v}}

\DeclareMathOperator{\re}{Re}
\DeclareMathOperator{\im}{Im}

\begin{document}
%
\title{Tight-and-Cheap Conic Relaxation for the\\Optimal Reactive Power Dispatch Problem}
%
%
%

\author{Christian~Bingane,~\IEEEmembership{Student~Member,~IEEE,}
        Miguel~F.~Anjos,~\IEEEmembership{Senior~Member,~IEEE}
        and~S\'ebastien~Le~Digabel 
\thanks{The second author is with the School of Mathematics, University of Edinburgh, 
Edinburgh EH9 3FD, Scotland, UK.}
\thanks{All authors are with the Department of Mathematics and Industrial Engineering, Polytechnique Montreal, Montreal, Quebec, Canada H3C~3A7; and the GERAD research center, Montreal, Quebec, Canada H3T~2A7.}
\thanks{E-mails: christian.bingane@polymtl.ca, anjos@stanfordalumni.org, sebastien.le-digabel@polymtl.ca.}
\thanks{This research was supported by the NSERC-Hydro-Quebec-Schneider Electric Industrial Research Chair.}
}

\maketitle

\begin{abstract}
The optimal reactive power dispatch (ORPD) problem is an alternating current optimal power flow (ACOPF) problem where discrete control devices for regulating the reactive power, such as shunt elements and tap changers, are considered. The ORPD problem is modelled as a mixed-integer nonlinear optimization problem and its complexity is increased compared to the ACOPF problem, which is highly nonconvex and generally hard to solve. Recently, convex relaxations of the ACOPF problem have attracted a significant interest since they can lead to global optimality. We propose a tight conic relaxation of the ORPD problem and show that a round-off technique applied with this relaxation leads to near-global optimal solutions with very small guaranteed optimality gaps, unlike with the nonconvex continuous relaxation. We report computational results on selected MATPOWER test cases with up to 3375 buses.
\end{abstract}


%
\IEEEpeerreviewmaketitle

\section*{Nomenclature}
\subsection{Notations}
\begin{abbrv}
\item[$\R$/$\C$] Set of real/complex numbers,
\item[$\Herm^n$] Set of $n\times n$ Hermitian matrices,
\item[$\jc$] Imaginary unit,
\item[$a$/$\cplx{a}$] Real/complex number,
\item[$\rv{a}$/$\cv{a}$] Real/complex vector,
\item[$A$/$\cplx{A}$]  Real/complex matrix.
\end{abbrv}
\subsection{Operators}
\begin{abbrv}
\item[$\re (\cdot)$/$\im (\cdot)$]  Real/imaginary part operator,
\item[$(\cdot)^*$] Conjugate operator,
\item[$\left|\cdot\right|$]  Magnitude or cardinality set operator,
\item[$\angle (\cdot)$]  Phase operator,
\item[$(\cdot)^H$]  Conjugate transpose operator.
\end{abbrv}
\subsection{Input data}
\begin{abbrv}
\item[$\Reseau = (\Node, \Branch)$] Power network,
\item[$\Node$] Set of buses,
\item[$\Shunt \subseteq \Node$] Set of buses $k$ where a shunt element is connected,
\item[$\Gen = \bigcup_{k\in\Node} \Gen_k$] Set of generators,
\item[$\Gen_k$] Set of generators connected to bus~$k$,
\item[$\Branch$] Set of branches,
\item[$\Tap \subseteq \Branch$] Set of branches with tap changers,
\item[$p_{Dk}$/$q_{Dk}$] Active/reactive power demand at bus~$k$,
\item[$g_k'$/$b_k'$] Conductance/susceptance of shunt element at bus~$k$,
\item[$c_{g2}, c_{g1}, c_{g0}$] Generation cost coefficients of generator~$g$,
\item[$\cplx{y}_\ell^{-1} = r_{\ell} + \jc x_{\ell}$] Series impedance of branch~$\ell$,
\item[$b'_\ell$] Total shunt susceptance of branch~$\ell$.
\end{abbrv}
\subsection{Variables}
\begin{abbrv}
\item[$p_{Gg}$/$q_{Gg}$] Active/reactive power generation by generator~$g$,
\item[$\cplx{v}_{k}$] Complex (phasor) voltage at bus~$k$,
\item[$u_{k}$] Shunt variable corresponding to the shunt element connected at bus~$k$,
\item[$p_{f\ell}$/$q_{f\ell}$] Active/reactive power flow injected along branch~$\ell$ by its \emph{from} end,
\item[$p_{t\ell}$/$q_{t\ell}$] Active/reactive power flow injected along branch~$\ell$ by its \emph{to} end,
\item[$t_{\ell}$] Turns ratio of tap changer~$\ell$.
\end{abbrv}

\section{Introduction}
\IEEEPARstart{T}{he} optimal power flow (OPF) problem, first formulated in~\cite{ref2}, consists in finding a network operating point that optimizes an objective function subject to power flow equations and other operational constraints~\cite{cain, frank1, frank2, capitanescu2016}. The continuous version with AC power flow equations, also called AC optimal power flow (ACOPF) problem, is nonconvex and NP-hard~\cite{verma2010,lehmann2016}. The optimal reactive power dispatch (ORPD) problem, also known as Volt/VAR optimization problem, can be seen as an ACOPF problem with discrete control devices for regulating the reactive power such as shunt elements and tap changers~\cite{yang2017,yang2017orpd}. Because of the presence of discrete variables, the ORPD problem is generally more difficult than the ACOPF.

A mixed-integer nonlinear program (MINLP) is an optimization problem which involves both continuous and integer variables and whose objective function and feasible set are described by nonlinear functions~\cite{bonami2012}. A MINLP is said to be \emph{convex} if its continuous relaxation, i.e. the problem obtained by dropping the integrality constraints, is a convex optimization problem; otherwise, it is said to be \emph{nonconvex}. MINLPs inherit difficulties from nonlinear programs (NLPs) and mixed-integer linear programs (MILPs) since they are a generalization of both classes. In fact, there exist simple cases of nonconvex MINLPs which are not only NP-hard, but even undecidable. See more details in~\cite{lee2011, burer2012}.

Recently, convex relaxations of the ACOPF problem, in particular second-order cone programming (SOCP)~\cite{jabr1}, semidefinite programming (SDP)~\cite{bai}, and quadratic convex (QC)~\cite{coffrin2016a} relaxations, have attracted a significant interest for several reasons. First, they can lead to global optimality; second, because they are relaxations, they provide a bound on the global optimal value of the ACOPF problem; and third, if one of these relaxations is infeasible, then the ACOPF problem is infeasible. Mixed-integer QC relaxations were proposed in~\cite{hijazi2014} and mixed-integer SOCP relaxations in~\cite{ding2016, wu2017, kocuk2017} for MINLPs in power systems. However, \cite{capitanescu2016} emphasized that convex relaxations of the OPF problem are aimed at complementing nonconvex (local) solvers with the valuable information about the quality of the solution obtained, rather than at replacing them.

Moreover, several applications of OPF are multi-period problems by nature due to factors such as changing market prices, ramping limits of generation units, and demand behavior~\cite{gemine2014}. Extending a convex relaxation from a single-period OPF problem to a multi-period one may jeopardize its exactness. This is discussed in~\cite{bingane2020} where a tight convex relaxation for the multi-period case is proposed.

On the other hand, some heuristic techniques discussed in~\cite{capitanescu2010} have been proposed for handling discrete variables in the OPF problem. One approach works as follows: first, solve the continuous relaxation of the OPF problem, i.e. treat the discrete variables as continuous; second, round-off solutions corresponding to discrete variables to their nearest discrete values; and third, fix discrete variables to these values and then solve the corresponding OPF subproblem. This approach is called the \emph{round-off technique} and remains the simplest to deal with discrete variables in the ORPD problem, although it may lead to poor suboptimal solutions or infeasible ones. Some deficiences were already pointed out in~\cite{tinney1988}: for instance, since the OPF problem is in general highly nonconvex, solving its continuous relaxation with a local optimizer in the first step may lead to a poor local solution which, after round-off in the second step, may lead to a very poor solution in the third step. More discussion about the round-off strategy can be found in~\cite{capitanescu2016, yang2017orpd}.

The main contribution of this paper is to show that a round-off technique used with a tight convex continuous relaxation may lead to near-global optimal solutions of the ORPD problem. This contribution is in two parts. First, we propose SDP-based relaxations of the ORPD problem with a new tight convex model of tap changer; second, a modified round-off technique where, in the first step, we solve a SDP-based relaxation instead of the nonconvex continuous relaxation of the ORPD problem. More details about semidefinite optimization can be found in~\cite{anjos2012}. Computational results show that the SDP-based relaxations of the ORPD problem are tight, and furthermore that applying the round-off technique with these relaxations leads to near-global solutions, even for large-scale instances. To the best of our knowledge, it is the first time that extensive computations with this approach are carried out for the ORPD problem with large-scale meshed networks.

The remainder of this paper is organized as follows. In Section~\ref{sec2:formulation}, we state the mathematical model of the ORPD problem. In Section~\ref{sec2:cvx}, we describe two SDP-based relaxations of the ORPD problem: a simple one already proposed in~\cite{lav1} and a new tighter one. In Section~\ref{sec2:algo}, we present our round-off technique, and we report computational results in Section~\ref{sec2:results}. 
Section~\ref{sec2:conclusion} concludes the paper.

\section{ORPD: Formulation}\label{sec2:formulation}
Consider a typical power network $\Reseau = (\Node, \Branch)$, where $\Node = \{1,2,\ldots,n\}$ and $\Branch \subseteq \Node \times \Node$ denote respectively the set of buses and the set of branches (transmission lines and tap changers). We denote by $\Shunt \subseteq \Node$ the set of buses $k$ where a shunt element is connected, and $\Tap \subseteq \Branch$ the set of branches with tap changers. Each branch~$\ell \in \Branch$ has a \emph{from} end  $k$ (on the \emph{tap side}) and a \emph{to} end $m$  as modeled in~\cite{matpower}. We denote $\ell = (k,m)$.  The ORPD problem is given as:
\begin{subequations}\label{eq2:orpd}
\begin{equation}\label{eq2:objfun}
\min f(\rv{u},\rv{t},\rv{p}_{G}, \rv{q}_{G},\rv{p}_{f},\rv{q}_{f},\rv{p}_{t},\rv{q}_{t}, \cv{v}) 
\end{equation}
over variables $\rv{u}\in \{0,1\}^{|\Shunt|}$, $\rv{t}\in \R^{|\Tap|}$, $\rv{p}_{G}, \rv{q}_{G} \in \R^{|\Gen|}$, $\rv{p}_{f},\rv{q}_{f},\rv{p}_{t},\rv{q}_{t} \in \R^{|\Branch|}$, and $\cv{v}\in~\C^{|\Node|}$, subject to
\begin{itemize}
\item Power balance equations:
\begin{align}
&\sum_{g\in \Gen_k} p_{Gg} - p_{Dk} - g_k'u_k\left|\cplx{v}_k\right|^2 \nonumber\\
& = \sum_{\ell=(k,m)\in\Branch} p_{f\ell} + \sum_{\ell=(m,k)\in\Branch} p_{t\ell}\;\forall k\in\Shunt,\label{eq2:kclpu}\\
&\sum_{g\in \Gen_k} q_{Gg} - q_{Dk} + b_k'u_k\left|\cplx{v}_k\right|^2 \nonumber\\
& = \sum_{\ell=(k,m)\in\Branch} q_{f\ell} + \sum_{\ell=(m,k)\in\Branch} q_{t\ell}\;\forall k\in\Shunt,\label{eq2:kclqu}\\
&\sum_{g\in \Gen_k} p_{Gg} - p_{Dk} \nonumber\\
& = \sum_{\ell=(k,m)\in\Branch} p_{f\ell} + \sum_{\ell=(m,k)\in\Branch} p_{t\ell}\;\forall k\in\Node\setminus\Shunt,\label{eq2:kclp}\\
&\sum_{g\in \Gen_k} q_{Gg} - q_{Dk} \nonumber\\
& = \sum_{\ell=(k,m)\in\Branch} q_{f\ell} + \sum_{\ell=(m,k)\in\Branch} q_{t\ell}\;\forall k\in\Node\setminus\Shunt,\label{eq2:kclq}
\end{align}
\item Branch flow equations:
\begin{align}
& \frac{\cplx{v}_k}{t_{\ell}} \left[\left(\jc \frac{b'_{\ell}}{2} + \cplx{y}_\ell \right) \frac{\cplx{v}_k}{t_{\ell}} - \cplx{y}_\ell \cplx{v}_m\right]^*\nonumber\\
&= p_{f\ell} + \jc q_{f\ell}\;\forall \ell = (k,m)\in\Tap, \label{eq2:sft}\\
& \cplx{v}_m \left[- \cplx{y}_\ell \frac{\cplx{v}_k}{t_{\ell}} + \left(\jc \frac{b'_{\ell}}{2} + \cplx{y}_\ell \right) \cplx{v}_m\right]^*\nonumber\\
&= p_{t\ell} + \jc q_{t\ell}\;\forall \ell = (k,m)\in\Tap, \label{eq2:stt}\\
& \cplx{v}_k \left[\left(\jc \frac{b'_{\ell}}{2} + \cplx{y}_\ell \right) \cplx{v}_k - \cplx{y}_\ell \cplx{v}_m\right]^*\nonumber\\
&= p_{f\ell} + \jc q_{f\ell}\;\forall \ell = (k,m)\in\Branch\setminus\Tap, \label{eq2:sf}\\
& \cplx{v}_m \left[- \cplx{y}_\ell \cplx{v}_k + \left(\jc \frac{b'_{\ell}}{2} + \cplx{y}_\ell \right) \cplx{v}_m\right]^*\nonumber\\
&= p_{t\ell} + \jc q_{t\ell}\;\forall \ell = (k,m)\in\Branch\setminus\Tap, \label{eq2:st}
\end{align}
\item Generator power capacities:
\begin{equation}\label{eq2:genlim}
\pl_{Gg}\le p_{Gg} \le \pu_{Gg},\, \ql_{Gg}\le q_{Gg} \le \qu_{Gg}\;\forall g\in\Gen,
\end{equation}
\item Line thermal limits:
\begin{equation}\label{eq2:linelim}
|p_{f\ell} + \jc q_{f\ell}| \le \overline{s}_{\ell},\, |p_{t\ell} + \jc q_{t\ell}| ~\le ~\overline{s}_{\ell}\;\forall \ell \in\Branch,
\end{equation}
\item Voltage magnitude limits:
\begin{equation}\label{eq2:magbus}
\vl_k\le \left|\cplx{v}_k\right| \le \vu_k\;\forall k\in\Node,
\end{equation}
\item Reference bus constraint:
\begin{equation}\label{eq2:slack}
\angle \cplx{v}_1 = 0,
\end{equation}
\item Shunt variable:
\begin{equation}\label{eq2:shunt}
 u_k \in \{0,1\}\;\forall k\in\Shunt,
\end{equation}
\item Tap ratio constraint:
\begin{equation}\label{eq2:tap}
t_\ell \in \{\lb{t}_\ell,\ldots,\ub{t}_\ell\}\;\forall \ell \in\Tap.
\end{equation}
\end{itemize}
\end{subequations}

Constraints \eqref{eq2:kclpu}--\eqref{eq2:st} derive from Kirchhoff's laws and represent power flow equations in the network. We assume $\vl_k > 0$ for all $k\in \Node$ in~\eqref{eq2:magbus} and $\lb{t}_\ell > 0$ for all $\ell \in \Tap$ in~\eqref{eq2:tap}. Constraint~\eqref{eq2:slack} specifies bus $k=1$ as the reference bus. In~\eqref{eq2:shunt}, we assume that the shunt connected to a bus~$k\in\Shunt$ has an on/off switch. For all $\ell \in\Tap$, the tap ratio $t_\ell$ in~\eqref{eq2:tap} is a discrete variable that typically takes on $2\ub{\eta} +1$ values $\{\hat{t}_{-\ub{\eta}},\ldots, \hat{t}_0, \ldots, \hat{t}_{\ub{\eta}}\}$, uniformly distributed around $\hat{t}_0$. In this paper, we assume without loss of generality that $\hat{t}_\ell$ varies between 0.9 and 1.1 pu with steps of 0.0125 pu (17 different settings) for all $\ell \in\Tap$.

The objective function $f(\rv{u},\rv{t},\rv{p}_{G}, \rv{q}_{G},\rv{p}_{f},\rv{q}_{f},\rv{p}_{t},\rv{q}_{t}, \cv{v})$ in~\eqref{eq2:objfun} may represent any objective function related to the ORPD problem: power loss, voltage deviation, number of control actions, generation cost, etc. Some of them are considered in~\cite{ke2018}. We note that power loss is a widely used objective function for the ORPD problem. However, according to~\cite{cain}, just minimizing loss is inconsistent with economic principles and may result in suboptimal dispatch while minimizing cost would be the correct objective function for economically dispatching resources and would inherently meet the objective of minimizing loss. Recently, \cite{yang2017orpd}~considered minimizing generation cost instead of power loss in their formulation of the ORPD problem. In this work, we consider both cost minimization and loss minimization in our computational results. We also assume that~\eqref{eq2:objfun} is a convex function, and then relaxations proposed in this paper remain valid for any convex objective function related to the ORPD problem.

\section{ORPD: Convexification}\label{sec2:cvx}
In this section, we describe two semidefinite relaxations of the ORPD problem: a simple one (SDR1) already proposed in~\cite{lav1} and a new tighter one (SDR2). SDR2 is obtained by combining SDR1 with a new tight convex model of tap changer and we show that SDR2 is then stronger than SDR1. Since both SDR1 and SDR2 can be expensive to solve for large-scale instances, we derive from them two relaxations that are cheaper to solve: TCR1 and TCR2 respectively.

\subsection{Semidefinite relaxation 1 (SDR1)}\label{subsec2:sdr1}
Let
\begin{subequations}
\begin{align}
\cplx{V} & := \cv{v}\cv{v}^H, \label{eq2:Vv}\\
\cplx{W}_{\{\ell\}} & :=
\begin{bmatrix}
\cplx{V}_{kk} & \cplx{W}_{k\ell} & \cplx{V}_{km}\\
\cplx{W}_{k\ell}^* & \cplx{W}_{\ell\ell} & \cplx{W}_{\ell m}\\
\cplx{V}_{km}^* & \cplx{W}_{\ell m}^* & \cplx{V}_{mm}
\end{bmatrix} \nonumber\\
& := \begin{bmatrix}
\cplx{v}_{k}\\
\cplx{w}_{\ell}\\
\cplx{v}_{m}
\end{bmatrix}
\begin{bmatrix}
\cplx{v}_{k}\\
\cplx{w}_{\ell}\\
\cplx{v}_{m}
\end{bmatrix}^H &\forall \ell = (k,m)\in\Tap, \label{eq2:Ww}\\
\cplx{w}_\ell & := \frac{\cplx{v}_k}{t_\ell} &\forall \ell = (k,m)\in\Tap. \label{eq2:wvt}
\end{align}
\end{subequations}
The ORPD problem~\eqref{eq2:orpd} can be reformulated as follows
\begin{subequations}\label{eq2:orpdVW}
\begin{align}
\text{minimize } &\eqref{eq2:objfun} \nonumber\\
\text{subject to } & \eqref{eq2:kclp}, \eqref{eq2:kclq}, \eqref{eq2:genlim}, \eqref{eq2:linelim}, \eqref{eq2:slack}, \eqref{eq2:shunt}, \eqref{eq2:tap}, \nonumber\\
& \eqref{eq2:Vv}, \eqref{eq2:Ww}, \eqref{eq2:wvt}, \nonumber\\
&\sum_{g\in \Gen_k} p_{Gg} - p_{Dk} - g_k'u_k\cplx{V}_{kk} \nonumber\\
& = \sum_{\ell=(k,m)\in\Branch} p_{f\ell} + \sum_{\ell=(m,k)\in\Branch} p_{t\ell}\;\forall k\in\Shunt, \label{eq2:kclpuV}\\
&\sum_{g\in \Gen_k} q_{Gg} - q_{Dk} + b_k'u_k\cplx{V}_{kk} \nonumber\\ & = \sum_{\ell=(k,m)\in\Branch} q_{f\ell} + \sum_{\ell=(m,k)\in\Branch} q_{t\ell}\;\forall k\in\Shunt,\label{eq2:kclquV}\\
&\left(-\jc \frac{b'_{\ell}}{2} + \cplx{y}_\ell^* \right) \cplx{W}_{\ell \ell} - \cplx{y}_\ell^* \cplx{W}_{\ell m} \nonumber\\
&= p_{f\ell} + \jc q_{f\ell}\;\forall \ell = (k,m) \in \Tap,\label{eq2:sfW}\\
&-\cplx{y}_\ell^* \cplx{W}_{\ell m}^* + \left(-\jc \frac{b'_{\ell}}{2} + \cplx{y}_\ell^* \right) \cplx{V}_{mm} \nonumber\\
&= p_{t\ell} + \jc q_{t\ell}\;\forall \ell = (k,m) \in \Tap,\label{eq2:stW}\\
&\left(-\jc \frac{b'_{\ell}}{2} + \cplx{y}_\ell^* \right) \cplx{V}_{kk} - \cplx{y}_\ell^* \cplx{V}_{km} \nonumber\\
&= p_{f\ell} + \jc q_{f\ell}\;\forall \ell = (k,m) \in \Branch \setminus \Tap, \label{eq2:sfV}\\
&-\cplx{y}_\ell^* \cplx{V}_{km}^* + \left(-\jc \frac{b'_{\ell}}{2} + \cplx{y}_\ell^* \right) \cplx{V}_{mm} \nonumber\\
&= p_{t\ell} + \jc q_{t\ell}\;\forall \ell = (k,m) \in \Branch \setminus \Tap, \label{eq2:stV}\\
&\vl_k^2\le \cplx{V}_{kk} \le \vu_k^2\;\forall k\in\Node.\label{eq2:Vkk}
\end{align}
\end{subequations}

If we define $\xi_k := u_k \cplx{V}_{kk} \in \{0,\cplx{V}_{kk}\}$ for all $k\in\Shunt$, a linear formulation of \eqref{eq2:kclpuV}--\eqref{eq2:kclquV} is given as
\begin{subequations}\label{eq2:uV1pq}
\begin{align}
&\sum_{g\in \Gen_k} p_{Gg} - p_{Dk} - g_k'\xi_k \nonumber\\
& = \sum_{\ell=(k,m)\in\Branch} p_{f\ell} + \sum_{\ell=(m,k)\in\Branch} p_{t\ell}\;\forall k\in\Shunt,\label{eq2:uV1p}\\
&\sum_{g\in \Gen_k} q_{Gg} - q_{Dk} + b_k'\xi_k \nonumber\\
& = \sum_{\ell=(k,m)\in\Branch} q_{f\ell} + \sum_{\ell=(m,k)\in\Branch} q_{t\ell}\;\forall k\in\Shunt,\label{eq2:uV1q}\\
& 0 \le \xi_k \le \cplx{V}_{kk}\;\forall k\in\Shunt.\label{eq2:uV1}
\end{align}
\end{subequations}
From \eqref{eq2:Ww} and \eqref{eq2:wvt}, we have
\begin{subequations}\label{eq2:WVt}
\begin{align}
\cplx{W}_{k\ell} &= \frac{\cplx{V}_{kk}}{t_\ell} &\forall \ell = (k,m)\in\Tap,\\
\cplx{W}_{\ell\ell} &= \frac{\cplx{V}_{kk}}{t_\ell^2} &\forall \ell = (k,m)\in\Tap,
\end{align}
\end{subequations}
which implies
\begin{subequations}\label{eq2:WVt1}
\begin{align}
\frac{\cplx{V}_{kk}}{\ub{t}_\ell} \le \cplx{W}_{k\ell} &\le \frac{\cplx{V}_{kk}}{\lb{t}_\ell} &\forall \ell = (k,m)\in\Tap,\\
\frac{\cplx{V}_{kk}}{\ub{t}_\ell^2} \le \cplx{W}_{\ell\ell} &\le \frac{\cplx{V}_{kk}}{\lb{t}_\ell^2} &\forall \ell = (k,m)\in\Tap,
\end{align}
\end{subequations}
since $0 < \lb{t}_\ell \le t_\ell \le \ub{t}_\ell$ for all $\ell \in \Tap$. Finally, we can show that $\cplx{V}$ in~\eqref{eq2:Vv} and, for all $\ell \in \Tap$, $\cplx{W}_{\{\ell\}}$ in~\eqref{eq2:Ww} are rank-one positive semidefinite matrices. The \emph{semidefinite relaxation 1} (SDR1) in Model~\ref{model2:sdr1} is obtained by dropping the rank constraints. For all $\ell = (k,m) \in \Tap$, since $\cplx{W}_{\{\ell\}} \succeq 0$ implies $\cplx{W}_{k\ell}^2 \le \cplx{V}_{kk} \cplx{W}_{\ell\ell}$, \eqref{eq2:WVt1}~can be rewritten as
\begin{subequations}\label{eq2:WVt1r}
\begin{align}
\cplx{W}_{k\ell} &\ge \frac{\cplx{V}_{kk}}{\ub{t}_\ell} &\forall \ell = (k,m)\in\Tap,\\
\cplx{W}_{\ell\ell} &\le \frac{\cplx{V}_{kk}}{\lb{t}_\ell^2} &\forall \ell = (k,m)\in\Tap.
\end{align}
\end{subequations}
\begin{algorithm}
\caption{Semidefinite relaxation 1 (SDR1)}
\label{model2:sdr1}
\begin{algorithmic}
\STATE Variables:
\begin{subequations}\label{eq2:var}
\begin{align}
\rv{\xi} &\in \R^{|\Shunt|},\\
\rv{p}_G, \rv{q}_G &\in \R^{|\Gen|},\\
\rv{p}_f, \rv{q}_f, \rv{p}_t, \rv{q}_t &\in\R^{|\Branch|},\\
\cplx{V} & \in\Herm^{|\Node|},\\
(\cplx{W}_{\ell\ell}, \cplx{W}_{k\ell},\cplx{W}_{\ell m}) & \in \R \times \R \times \C & \forall \ell = (k,m) \in \Tap.
\end{align}
\end{subequations}
\STATE Minimize:~\eqref{eq2:objfun}
\STATE Subject to: \eqref{eq2:kclp}, \eqref{eq2:kclq}, \eqref{eq2:genlim}, \eqref{eq2:linelim}, \eqref{eq2:sfW}--\eqref{eq2:Vkk}, \eqref{eq2:uV1pq}, \eqref{eq2:WVt1r}, $\cplx{V}\succeq 0$, and $\cplx{W}_{\{\ell\}} \succeq 0$ $\forall \ell \in\Tap$.
\end{algorithmic}
\end{algorithm}

We point out that the variables $u_k$ and $t_\ell$ are eliminated in Model~\ref{model2:sdr1} for all $k\in\Shunt$ and for all $\ell \in\Tap$. This approach was already proposed in~\cite{lav1} and was recently used in~\cite{robbins2016}. Once the solution of Model~\ref{model2:sdr1} is obtained, the optimal value~$\hat{u}_k$ of the shunt element connected to a bus~$k\in\Shunt$ and the optimal tap ratio $\hat{t}_\ell$ of a transformer $\ell = (k,m) \in \Tap$ can be determined respectively as follows:
\begin{subequations}\label{eq2:ut}
\begin{align}
\hat{u}_k = \frac{\xi_k}{\cplx{V}_{kk}} &\in [0,1],\\
\hat{t}_\ell = \sqrt{\frac{\cplx{V}_{kk}}{\cplx{W}_{\ell\ell}}} &\in [\lb{t}_\ell, \ub{t}_\ell].
\end{align}
\end{subequations}

If the optimal solutions $\hat{\cplx{V}}$ and $\hat{\cplx{W}}_{\{\ell\}}$ for all $\ell \in\Tap$ are rank-one matrices, $\hat{u}_k \in \{0,1\}$ for all $k\in\Shunt$, and $\hat{t}_\ell \in \{\lb{t}_\ell,\ldots,\ub{t}_\ell\}$ for all $\ell \in\Tap$, then there exists a complex vector $\hat{\cv{v}}$ that is a global optimal solution of~\eqref{eq2:orpd}. We say that the relaxation SDR1 is exact.

\begin{subequations}\label{eq2:tcr}
We also consider a cheaper relaxation, called ``\emph{tight-and-cheap relaxation 1}'', given in Model~\ref{model2:tcr1} and obtained as follows. We replace~$\cplx{V} \succeq 0$ and~$\cplx{W}_{\{\ell\}} \succeq 0$ for all $\ell \in \Tap$ in Model~\ref{model2:sdr1} by
\begin{align}
\begin{bmatrix}
1 & \cplx{v}_{k}^* & \cplx{v}_{m}^*\\
\cplx{v}_{k} &\cplx{V}_{kk} & \cplx{V}_{km}\\
\cplx{v}_{m} & \cplx{V}_{km}^* & \cplx{V}_{mm}
\end{bmatrix} &\succeq 0 &\forall \ell = (k,m)\in\Branch \setminus \Tap, \label{eq2:Vvsdp}\\
\begin{bmatrix}
1 & \cv{w}_{\{\ell\}}^H \\
\cv{w}_{\{\ell\}} &\cplx{W}_{\{\ell\}}
\end{bmatrix} &\succeq 0 &\forall \ell = (k,m)\in\Tap, \label{eq2:Wwsdp}
\end{align}
where $\cv{w}_{\{\ell\}}^H =(\cplx{v}_k^*,\cplx{w}_{\ell}^*,\cplx{v}_m^*)$ for all $\ell = (k,m)\in\Tap$, and we add the following constraints
\begin{align}\label{eq2:rlt}
\re (\cplx{v}_1) &\ge \frac{\cplx{V}_{11} + \vl_1 \vu_1}{\vl_1 + \vu_1},\\
\im (\cplx{v}_1) &= 0,
\end{align}
corresponding to the reference bus~$k=1$.
\end{subequations}
\begin{algorithm}
\caption{Tight-and-cheap relaxation 1 (TCR1)}
\label{model2:tcr1}
\begin{algorithmic}
\begin{subequations}
\STATE Variables:~\eqref{eq2:var}, $\cv{v} \in \C^{|\Node|}$, $\cv{w} \in \C^{|\Tap|}$.
\STATE Minimize:~\eqref{eq2:objfun}
\STATE Subject to: \eqref{eq2:kclp}, \eqref{eq2:kclq}, \eqref{eq2:genlim}, \eqref{eq2:linelim}, \eqref{eq2:sfW}--\eqref{eq2:Vkk}, \eqref{eq2:uV1pq}, \eqref{eq2:WVt1}, \eqref{eq2:tcr}.
\end{subequations}
\end{algorithmic}
\end{algorithm}

The tight-and-cheap relaxation (TCR) was first proposed in~\cite{bingane2018} for the ACOPF problem. It was shown in~\cite{bingane2018} that TCR is stronger than the standard SOCP relaxation and nearly as tight as the standard SDP relaxation. Moreover, computation experiments on standard test cases with up to 6515 buses showed that solving TCR  for large-scale instances is much less expensive than solving the chordal relaxation, a SDP relaxation technique that exploits the sparsity of power networks.

\subsection{Semidefinite relaxation 2 (SDR2)}
For all $\ell = (k,m) \in \Tap$, variables $\cplx{W}_{k\ell}$ and $\cplx{W}_{\ell\ell}$ in~\eqref{eq2:WVt} are respectively described by constraints of the form $z_n = x/y^n$, $n=1,2$. Consider the set
\[
\begin{aligned}
\mathcal{S}_1 = \{(x,y,z_1,z_2)\in \R^4 \colon &\lb{x} \le x \le \ub{x}, \lb{y} \le y \le \ub{y},\\
& z_1 = x/y, z_2 = x/y^2\},
\end{aligned}
\]
where $0 < \lb{x} < \ub{x}$ and $ 0 < \lb{y} < \ub{y}$. We can show that $\mathcal{S}_1$ is equivalent to
\[
\mathcal{S}_2 = \{(x,z_1,z_2)\in \R^3 \colon z_2 = z_1^2/x, (x,z_1) \in \Omega\},
\]
where $\Omega = \{(x,z_1)\in \R^2 \colon \lb{x} \le x \le \ub{x}, x/\ub{y} \le z_1\le x/\lb{y}\}$ is the convex quadrilateral with vertices $(\lb{x}, \lb{x}/\ub{y})$, $(\ub{x}, \ub{x}/\ub{y})$, $(\ub{x}, \ub{x}/\lb{y})$, and $(\lb{x}, \lb{x}/\lb{y})$.

In Section~\ref{subsec2:sdr1}, the proposed convex set that contains~$\mathcal{S}_2$ is
\[
\begin{aligned}
\mathcal{S}_3 = \{(x,z_1,z_2)\in \R^3 \colon &\lb{x} \le x \le \ub{x},\\
& z_2 \ge z_1^2/x, z_1 \ge x/\ub{y}, z_2 \le x/\lb{y}^2\}.
\end{aligned}
\]
In this section, we propose a tighter convex set containing~$\mathcal{S}_2$.
\begin{lemma}\label{lemma2:x2/yconv}
Let $f(x,z_1) = z_1^2/x$ defined over $\R_+ \times \R$. Then $f$ is convex.
\end{lemma}
\begin{IEEEproof}
For all $(x,z_1) \in \R_+ \times \R$,
\[
\nabla^2 f(x,z_1) = \frac{2}{x^3}
\begin{bmatrix}
z_1^2 & -xz_1\\
-xz_1 & x^2
\end{bmatrix} = \frac{2}{x^3}
\begin{bmatrix}
z_1\\
-x
\end{bmatrix} 
\begin{bmatrix}
z_1\\
-x
\end{bmatrix}^T \succeq 0.
\]
Then $f$ is convex.
\end{IEEEproof}
\begin{lemma}\label{lemma2:x2/yconc}
Let $(x,z_1) \in \Omega$. If $z_2 = z_1^2/x$, then $x + \lb{y} \ub{y} z_2 \le (\lb{y} + \ub{y})z_1$.
\end{lemma}
\begin{IEEEproof}
Let $(x,z_1) \in \Omega$ and $z_2 = z_1^2/x$. We have $x/\ub{y} \le z_1 \le x/\lb{y}$ with $0 < \lb{y} < \ub{y}$. Therefore, $(x - \lb{y}z_1)(x -\ub{y}z_1) \le 0$. On the other hand,
\[
\begin{aligned}
(x - \lb{y}z_1)(x -\ub{y}z_1) &= x^2 - (\lb{y} + \ub{y})xz_1 + \lb{y} \ub{y} z_1^2\\
&= x[x - (\lb{y} + \ub{y})z_1 + \lb{y} \ub{y}z_2] \le 0.
\end{aligned}
\]
Since $\lb{x} > 0$, it follows that $x - (\lb{y} + \ub{y})z_1 + \lb{y} \ub{y}z_2 \le 0$.
\end{IEEEproof}

\begin{proposition}\label{prop2:concenv}
Let $f(x,z_1) = z_1^2/x$ defined over~$\Omega$. The concave envelope or the tightest concave overestimator of~$f$ on~$\Omega$ is $h(x,z_1) = \frac{1}{\lb{y}\ub{y}}[(\lb{y} + \ub{y})z_1 - x]$.
\end{proposition}
\begin{IEEEproof}
The function~$h$ is affine and, from Lemma~\ref{lemma2:x2/yconc}, overestimates~$f$ on~$\Omega$. Then $h$ is a concave overestimator of~$f$ on~$\Omega$. We note that $h(x,z_1) = f(x,z_1)$ for all $(x,z_1) \in \Omega$ such that $z_1 = x/\ub{y}$ or $z_1 = x/\lb{y}$.

Suppose that there exists a concave overestimator~$\tilde{h}$ of~$f$ on~$\Omega$ such that $f(x,z_1) \le \tilde{h}(x,z_1) \le h(x,z_1)$ for all $(x,z_1) \in \Omega$ and $\tilde{h}(\tilde{x},\tilde{z}_1) < h(\tilde{x},\tilde{z}_1)$ for some $(\tilde{x},\tilde{z}_1) \in \Omega$. Since~$\Omega$ is the convex quadrilateral with vertices $(\lb{x}, \lb{x}/\ub{y})$, $(\ub{x}, \ub{x}/\ub{y})$, $(\ub{x}, \ub{x}/\lb{y})$, and $(\lb{x}, \lb{x}/\lb{y})$, there exist nonnegative scalars $\alpha_1$, $\alpha_2$, $\alpha_3$, and $\alpha_4$ such that $\alpha_1 + \alpha_2 + \alpha_3 + \alpha_4 = 1$ and $(\tilde{x},\tilde{z}_1) = \alpha_1 (\lb{x}, \lb{x}/\ub{y})+ \alpha_2 (\ub{x}, \ub{x}/\ub{y}) + \alpha_3 (\ub{x}, \ub{x}/\lb{y}) + \alpha_4 (\lb{x}, \lb{x}/\lb{y})$. Therefore,
\[
\begin{aligned}
\tilde{h}(\tilde{x},\tilde{z}_1) &= \tilde{h}(\alpha_1 (\lb{x}, \lb{x}/\ub{y})+ \alpha_2 (\ub{x}, \ub{x}/\ub{y})\\
&\phantom{=} + \alpha_3 (\ub{x}, \ub{x}/\lb{y}) + \alpha_4 (\lb{x}, \lb{x}/\lb{y}))\\
& \ge \alpha_1 \tilde{h}(\lb{x}, \lb{x}/\ub{y})+ \alpha_2 \tilde{h}(\ub{x}, \ub{x}/\ub{y})\\
&\phantom{=} + \alpha_3 \tilde{h}(\ub{x}, \ub{x}/\lb{y}) + \alpha_4 \tilde{h}(\lb{x}, \lb{x}/\lb{y})\\
& = \alpha_1 f(\lb{x}, \lb{x}/\ub{y})+ \alpha_2 f(\ub{x}, \ub{x}/\ub{y})\\
&\phantom{=} + \alpha_3 f(\ub{x}, \ub{x}/\lb{y}) + \alpha_4 f(\lb{x}, \lb{x}/\lb{y})\\
& = \alpha_1 h(\lb{x}, \lb{x}/\ub{y})+ \alpha_2 h(\ub{x}, \ub{x}/\ub{y})\\
&\phantom{=} + \alpha_3 h(\ub{x}, \ub{x}/\lb{y}) + \alpha_4 h(\lb{x}, \lb{x}/\lb{y})\\
& = h(\alpha_1 (\lb{x}, \lb{x}/\ub{y})+ \alpha_2 (\ub{x}, \ub{x}/\ub{y})\\
&\phantom{=} + \alpha_3 (\ub{x}, \ub{x}/\lb{y}) + \alpha_4 (\lb{x}, \lb{x}/\lb{y}))\\
& = h(\tilde{x},\tilde{z}_1).
\end{aligned}
\]
The inequality follows from the definition of a concave function; the two subsequent equalities from the fact that~$h$ agrees with~$f$ at the vertices of~$\Omega$ and $f(x,z_1) \le \tilde{h}(x,z_1) \le h(x,z_1)$ for all $(x,z_1) \in \Omega$; and the last equality from the fact that~$h$ is affine. The relation obtained $\tilde{h}(\tilde{x},\tilde{z}_1) \ge h(\tilde{x},\tilde{z}_1)$ contradicts $\tilde{h}(\tilde{x},\tilde{z}_1) < h(\tilde{x},\tilde{z}_1)$. So $h$ is the concave envelope of~$f$ on~$\Omega$.
\end{IEEEproof}

A similar proof can be found in~\cite{benson2004} on convex and concave envelopes of functions defined on convex quadrilaterals.

\begin{proposition}\label{prop2:convhull}
The convex hull of $\mathcal{S}_2 = \{(x,z_1,z_2)\in \R^3 \colon \lb{x} \le x \le \ub{x}, x/\ub{y} \le z_1 \le x/\lb{y}, z_2 = z_1^2/x\}$, where $0 < \lb{x} < \ub{x}$, $ 0 < \lb{y} < \ub{y}$, is $\ub{\mathcal{S}}_2 = \{(x,z_1,z_2)\in~\R^3 \colon \lb{x} \le x \le \ub{x},  z_1^2 \le xz_2, x + \lb{y} \ub{y} z_2 \le (\lb{y} + \ub{y})z_1\}$.
\end{proposition}
\begin{IEEEproof}
The result follows from Lemma~\ref{lemma2:x2/yconv} and Proposition~\ref{prop2:concenv}.
\end{IEEEproof}

Applying Proposition~\ref{prop2:convhull} to~\eqref{eq2:WVt}, we replace~\eqref{eq2:WVt1r} in Model~\ref{model2:sdr1} by
\begin{equation}\label{eq2:WVt2}
\cplx{V}_{kk} + \lb{t}_\ell \ub{t}_\ell \cplx{W}_{\ell\ell} \le (\lb{t}_\ell + \ub{t}_\ell) \cplx{W}_{k\ell}\; \forall \ell = (k,m) \in \Tap,
\end{equation}
and define a new formulation of SDP relaxation in Model~\ref{model2:sdr2}. We will refer to this relaxation as ``\emph{semidefinite relaxation~2}'' (SDR2). The ``\emph{tight-and-cheap relaxation~2}'' is given in Model~\ref{model2:tcr2}. From Proposition~\ref{prop2:convhull}, $\ub{\mathcal{S}}_2 \subset \mathcal{S}_3$ and then SDR2 (respectively TCR2) is stronger than SDR1 (TCR1).

\begin{algorithm}
\caption{Semidefinite relaxation 2 (SDR2)}
\label{model2:sdr2}
\begin{algorithmic}
\begin{subequations}
\STATE Variables:~\eqref{eq2:var}.
\STATE Minimize:~\eqref{eq2:objfun}
\STATE Subject to: \eqref{eq2:kclp}, \eqref{eq2:kclq}, \eqref{eq2:genlim}, \eqref{eq2:linelim}, \eqref{eq2:sfW}--\eqref{eq2:Vkk}, \eqref{eq2:uV1pq}, \eqref{eq2:WVt2}, $\cplx{V}\succeq 0$, and $\cplx{W}_{\{\ell\}} \succeq 0$ $\forall \ell \in\Tap$.
\end{subequations}
\end{algorithmic}
\end{algorithm}

\begin{algorithm}
\caption{Tight-and-cheap relaxation 2 (TCR2)}
\label{model2:tcr2}
\begin{algorithmic}
\begin{subequations}
\STATE Variables:~\eqref{eq2:var}, $\cv{v} \in \C^{|\Node|}$, $\cv{w} \in \C^{|\Tap|}$.
\STATE Minimize:~\eqref{eq2:objfun}
\STATE Subject to: \eqref{eq2:kclp}, \eqref{eq2:kclq}, \eqref{eq2:genlim}, \eqref{eq2:linelim}, \eqref{eq2:sfW}--\eqref{eq2:Vkk}, \eqref{eq2:uV1pq}, \eqref{eq2:WVt2}, \eqref{eq2:tcr}.
\end{subequations}
\end{algorithmic}
\end{algorithm}

\section{ORPD: Solution approach}\label{sec2:algo}
Following~\cite{burer2012}, a mixed-integer nonlinear program (MINLP) is an optimization problem of the form
\begin{subequations}\label{eq2:minlp}
\begin{align}
\text{$\hat{\upsilon}=$ minimize } & f_0 (\rv{x},\rv{y})\\
\text{subject to } & f_i (\rv{x},\rv{y}) \le 0, & i=1,2,\ldots,m,\\
& \rv{x}\in\R^{n_1}, \rv{y}\in\Z^{n_2}.
\end{align}
\end{subequations}

By fixing all integer variables, i.e. $\rv{y} = \tilde{\rv{y}}\in\Z^{n_2}$, we obtain the following \emph{subproblem} (SP) of~\eqref{eq2:minlp}
\begin{subequations}\label{eq2:sp}
\begin{align}
\text{$\hat{\upsilon}_{SP}=$ minimize } & f_0 (\rv{x},\tilde{\rv{y}})\\
\text{subject to } & f_i (\rv{x},\tilde{\rv{y}}) \le 0, & i=1,2,\ldots,m,\\
& \rv{x}\in\R^{n_1}.
\end{align}
\end{subequations}
Then we can rewrite the MINLP~\eqref{eq2:minlp} as
\[
\hat{\upsilon}= \min_{\tilde{\rv{y}}\in\Z^{n_2}} \left\{\min_{\rv{x}\in\R^{n_1}} \{f_0 (\rv{x},\tilde{\rv{y}}) \colon f_i (\rv{x},\tilde{\rv{y}}) \le 0, i=1,2,\ldots,m \}\right\}
\]
If~\eqref{eq2:sp} is feasible, its global optimal value $\hat{\upsilon}_{SP}$ provides an upper bound for~$\hat{\upsilon}$ and the optimal solution~$(\hat{\rv{x}}_{SP}, \tilde{\rv{y}})$ of~\eqref{eq2:sp} is called a \emph{suboptimal} solution of MINLP \eqref{eq2:minlp}. The subproblem~\eqref{eq2:sp} may be nonconvex and solving it to global optimality may be hard. In this case, any local optimal value~$\ub{\upsilon}_{SP}$ of~\eqref{eq2:sp} can be considered as an upper bound of~$\hat{\upsilon}$.

Consider a relaxation~(R) of~\eqref{eq2:minlp} of the form
\begin{subequations}\label{eq2:r}
\begin{align}
\text{$\hat{\upsilon}_{R}=$ minimize } & \lb{f}_0 (\rv{x},\rv{y})\\
\text{subject to } & (\rv{x},\rv{y})\in\ub{\Omega} \subseteq \R^{n_1} \times \R^{n_2},
\end{align}
\end{subequations}
where $\ub{\Omega}$ contains the feasible set of~\eqref{eq2:minlp} and $\lb{f}_0$ underestimates $f_0$ for all feasible solutions of~\eqref{eq2:minlp}. The global optimal value~$\hat{\upsilon}_{R}$ is a lower bound of~$\hat{\upsilon}$. Since it is necessary to solve~\eqref{eq2:r} to global optimality to obtain a lower bound on~$\hat{\upsilon}$, it is more advantageous for the relaxation~\eqref{eq2:r} to be convex.

Let~$(\hat{\rv{x}}_R, \hat{\rv{y}}_R)$ and~$\hat{\upsilon}_{R}$ be respectively the global optimal solution and the global optimal value of a convex relaxation~(R) of a MINLP. Now, consider the subproblem~(SP) associated to~$\tilde{\rv{y}}_R$, the closest integer solution to~$\hat{\rv{y}}_R$. We define the optimality gap of the relaxation (R) by $100(1-\hat{\upsilon}_{R}/\ub{\upsilon}_{SP})$, where~$\ub{\upsilon}_{SP}$ is a local optimal value of the subproblem (SP). If the optimality gap is close to zero, we say that the relaxation is \emph{tight} and the suboptimal solution~$(\hat{\rv{x}}_{SP}, \tilde{\rv{y}}_R)$ is \emph{near-global optimal} for the MINLP.

This leads to the following approach to solve the ORPD problem~\eqref{eq2:orpd}:
\begin{enumerate}
\item Solve SDR1, TCR1, SDR2 or TCR2 and find corresponding shunt solutions $\hat{\rv{u}} \in [0,1]^{|\Shunt|}$ and tap ratios $\hat{\rv{t}} \in \prod_{\ell \in \Tap} [\lb{t}_\ell, \ub{t}_\ell]$ with formulas~\eqref{eq2:ut},
\item Round-off $\rv{u}$ and $\rv{t}$ to their respective nearest discrete values $\tilde{\rv{u}} \in \{0,1\}^{|\Shunt|}$ and $\tilde{\rv{t}} \in \prod_{\ell \in \Tap} \{\lb{t}_\ell,\ldots, \ub{t}_\ell\}$,
\item Fix $\rv{u} = \tilde{\rv{u}}$ and $\rv{t} = \tilde{\rv{t}}$ and solve the ACOPF subproblem with a nonlinear (local) solver.
\end{enumerate}

\section{Computational results}\label{sec2:results}
In this section, we evaluate the accuracy and the computational efficiency of SDR1 and TCR1 as compared to SDR2 and TCR2.

In order to increase the computational speed of solving SDR1 and SDR2, we exploited the sparsity of a power network by replacing the SDP constraint~$\cplx{V} \succeq 0$ by small SDP constraints defined on a chordal extension of the power network. More details about exploiting the sparsity of power networks in SDP relaxations of the OPF problem can be found in~\cite{bingane2018, jabr4, low1, andersen2014, madani2014}.

Let us interpret the network $\Reseau = (\Node, \Branch)$ as a connected, simple and undirected graph $\G = (\Node,\E)$ where $\Node = \{1,\ldots,n\}$ represents the set of vertices and $\E =\{\{k,m\}: (k,m) \text{ or } (m,k) \in \Branch\}$, the set of edges. It was proved in~\cite{grone} that the SDP constraint $\cplx{V}\succeq 0$ in SDR1 or SDR2 is equivalent to $\cplx{V}_\K \succeq 0$ for every maximal clique~$\K$ of a chordal extension~$\G'$ of~$\G$. $\cplx{V}_\K$ is the submatrix of $\cplx{V}$ in which the set of row indices that remain and the set of column indices that remain are both~$\K$. To compute the maximal cliques of a chordal extension of~$\G$, we used the same algorithm as in~\cite{bingane2018}.

We tested the models~\ref{model2:sdr1},~\ref{model2:tcr1},~\ref{model2:sdr2},~\ref{model2:tcr2}, on standard test cases available from MATPOWER~\cite{matpower, birch, josz2016ac}. We assigned a shunt variable~$u_k$ to each shunt element connected to a bus~$k\in\Shunt$ and we assumed that the tap ratio $t_\ell$ of each transformer $\ell \in \Tap$ varies from 0.9 to 1.1~pu by steps of 0.0125~pu. Table~\ref{table2:data} lists the test cases along with the number of shunt elements~$|\Shunt|$ and the number of transformers~$|\Tap|$.

We solved SDR1, TCR1, SDR2 and TCR2 in MATLAB using CVX~2.1~\cite{cvx2} with the solver MOSEK~8.0.0.60 and \texttt{default precision} (tolerance $\epsilon = 1.49 \times 10^{-8}$). All relaxations were implemented as a MATLAB package, which is available on GitHub~\cite{conicopf}. It requires that MATPOWER and CVX be installed and that the input instance be provided in MATPOWER format. MOSEK numerically failed to solve SDR1 for \texttt{case3120sp}.

All ACOPF subproblems were solved with the MATPOWER-solver MIPS. When MIPS numerically failed to solve subproblems for some test cases, marked with~``*'' in Table~\ref{table2:cost1} and Table~\ref{table2:loss1}, we then used the solver FMINCON. Among these subproblems, FMINCON did not converge to a feasible solution after 1000 iterations for \texttt{case3012wp} (SDR1), \texttt{case3120sp} (TCR1) in cost minimization, and \texttt{case3120sp} (SDR1, TCR1), \texttt{case3375wp} (SDR1, TCR1) in loss minimization.

All the computations were carried out on an \texttt{Intel Core i7-6700 CPU @ 3.40 GHz} computing platform.

We report results for two different objective functions: the generation cost $\sum_{g\in\Gen} c_{g2} p_{Gg}^2 + c_{g1} p_{Gg} + c_{g0}$~[\$/h] and the active power losses $\sum_{g\in\Gen} p_{Gg}$~[MW]. Both objective functions of test cases from~\cite{josz2016ac} are the same. We denote $\lb{\upsilon}$ the best lower bound which is the maximum value among $\hat{\upsilon}_{SDR1}$, $\hat{\upsilon}_{SDR2}$, $\hat{\upsilon}_{TCR1}$, $\hat{\upsilon}_{TCR2}$, respective optimal values of SDR1, TCR1, SDR2, and TCR2. The normalized~$\hat{\upsilon}_R$ of a relaxation is measured as $\hat{\upsilon}_R/\lb{\upsilon} \le 1$, where $\hat{\upsilon}_R$ is the relaxation's optimal value. Table~\ref{table2:cost1} and Table~\ref{table2:loss1} summarize the normalized optimal values of relaxations. The results support the following key points:
\begin{enumerate}
\item Among all relaxations, SDR2 is the strongest.
\item SDR1 (respectively SDR2) is stronger than TCR1 (respectively TCR2).
\item SDR1 and TCR2 are comparable.
\end{enumerate}

For a relaxation, the local optimal value of the ACOPF subproblem obtained from the relaxation's optimal solution after rounding-off the discrete variables $\rv{u}$ and $\rv{t}$ is denoted~$\ub{\upsilon}_{SP}$. The normalized~$\ub{\upsilon}_{SP}$ in Table~\ref{table2:cost1} and Table~\ref{table2:loss1} is the value $\ub{\upsilon}_{SP}/\ub{\upsilon} \ge 1$, where $\ub{\upsilon}$ the best upper bound, i.e. the minimum value among all $\ub{\upsilon}_{SP}$. In Table~\ref{table2:cost1} and Table~\ref{table2:loss1}, we observe that:
\begin{enumerate}
\item In general, applying the round-off technique with SDR2 (respectively TCR2) provide tighter upper bounds than with SDR1 (respectively TCR1) for large-scale instances.
\item Applying the round-off technique with TCR2 is comparable to using SDR1.
\end{enumerate}

Optimality gaps of the four relaxations are given in Table~\ref{table2:cost2} and Table~\ref{table2:loss2}. The optimality gap of a relaxation is measured as $100(1-\hat{\upsilon}_R/\ub{\upsilon}_{SP})$, where $\hat{\upsilon}_R$ is the relaxation's optimal value and $\ub{\upsilon}_{SP}$ is a local optimal value of the corresponding subproblem. With the optimality gap, we can guarantee the near-global optimality of the suboptimal solution obtained with MIPS. Results in Table~\ref{table2:cost2} and Table~\ref{table2:loss2} show that:
\begin{enumerate}
\item Suboptimal solutions obtained from SDR2's optimal solutions have the lowest guaranteed optimality gaps. We can say that they are near-global optimal for all but one test case: \texttt{case\_ACTIV\_SG\_500} (cost minimization).
\item In general, suboptimal solutions obtained from SDR1's and TCR2's optimal solutions have almost the same guaranteed optimality gaps.
\item Suboptimal solutions obtained from TCR1's optimal solutions have slightly larger optimality gaps.
\end{enumerate}

The computation times required to solve SDR1, TCR1, SDR2, and TCR2 as reported by MOSEK are also shown in Table~\ref{table2:cost2} and Table~\ref{table2:loss2}. We see that
\begin{enumerate}
\item Solving SDR1 (respectively TCR1) is as expensive as solving SDR2 (respectively TCR2) in general.
\item Solving SDR1 (respectively SDR2) is much more expensive than solving TCR1 (respectively TCR2) for large-scale instances. The TCRs are on average 30 times faster than the SDRs.
\end{enumerate}

Overall, we can see that TCR2 offers an interesting trade-off between the optimality gap and the computation time.

\section{Conclusion}\label{sec2:conclusion}
We proposed a tight SDP relaxation (SDR2) for the ORPD problem. This formulation is based on a standard SDP relaxation (SDR1) combined with a tight convex tap changer model. Experiments on selected MATPOWER instances with up to 3375 buses show that SDR2 is stronger than SDR1, and computationally comparable. From both SDR1 and SDR2, we derived tight-and-cheap relaxations TCR1 and TCR2, respectively.

A round-off technique based on a SDP relaxation of the ORPD problem instead of the nonconvex continuous relaxation was also proposed. Computational results show that applying the round-off technique with SDR2 and TCR2 provides suboptimal solutions which are near-global optimal. Both provide almost the same guaranteed optimality gaps, but TCR2 is computationally much less expensive for large-scale instances. In summary, the proposed TCR-based approach provides the best trade-off between optimality gap and computation time compared to the SDR-based approach.

\begin{table}[!t]
\footnotesize
\centering
\caption{Dimensions of test instances}\label{table2:data}
\begin{tabular}{@{}lrrrr@{}}
\toprule
Test case & $|\Node|$ & $|\Branch|$ & $|\Shunt|$ & $|\Tap|$ \\
\midrule
\multicolumn{5}{c}{\it Small-scale instances}\\
\midrule
{\tt case14}	&	14	&	20	&	1	&	3	\\
{\tt case24\_ieee\_rts}	&	24	&	38	&	1	&	5	\\
{\tt case30}	&	30	&	41	&	2	&	0	\\
{\tt case\_ieee30}	&	30	&	41	&	2	&	4	\\
{\tt case39}	&	39	&	46	&	0	&	12	\\
{\tt case57}	&	57	&	80	&	3	&	17	\\
{\tt case89pegase}	&	89	&	210	&	44	&	32	\\
\midrule
\multicolumn{5}{c}{\it Medium-scale instances}\\
\midrule
{\tt case118}	&	118	&	186	&	14	&	9	\\
{\tt case\_ACTIV\_SG\_200}	&	200	&	245	&	4	&	66	\\
{\tt case\_illinois200}	&	200	&	245	&	4	&	66	\\
{\tt case300}	&	300	&	411	&	29	&	107	\\
{\tt case\_ACTIV\_SG\_500}	&	500	&	597	&	15	&	131	\\
\midrule
\multicolumn{5}{c}{\it Large-scale instances}\\
\midrule
{\tt case1354pegase}	&	1 354	&	1 991	&	1 082	&	234	\\
{\tt case2383wp}	&	2 383	&	2 896	&	0	&	164	\\
{\tt case2736sp}	&	2 736	&	3 269	&	1	&	168	\\
{\tt case2737sop}	&	2 737	&	3 269	&	5	&	169	\\
{\tt case2746wop}	&	2 746	&	3 307	&	6	&	171	\\
{\tt case2746wp}	&	2 746	&	3 279	&	0	&	171	\\
{\tt case2869pegase}	&	2 869	&	4 582	&	2 197	&	493	\\
{\tt case3012wp}	&	3 012	&	3 572	&	9	&	201	\\
{\tt case3120sp}	&	3 120	&	3 693	&	9	&	206	\\
{\tt case3375wp}	&	3 374	&	4 161	&	9	&	381	\\
\bottomrule
\end{tabular}
\end{table}

\begin{table*}[!t]
\footnotesize
\centering
\caption{Cost minimization: Normalized optimal values}\label{table2:cost1}
\begin{tabular}{@{}lrr|rrrr|rrrr@{}}
\toprule
Test case & $\lb{\upsilon}$ [\$/h] & $\ub{\upsilon}$ [\$/h] & \multicolumn{4}{c|}{Normalized $\hat{\upsilon}_R$} & \multicolumn{4}{c}{Normalized $\ub{\upsilon}_{SP}$} \\
\cmidrule(r){4-7} \cmidrule(l){8-11} &&& SDR1 & SDR2 & TCR1 & TCR2 & SDR1 & SDR2 &  TCR1 & TCR2 \\
\midrule
\multicolumn{11}{c}{\it Small-scale instances}\\
\midrule
{\tt case14}	&	8 078.62	&	8 078.75	&	1.0000	&	1.0000	&	1.0000	&	1.0000	&	1.0000	&	1.0000	&	1.0000	&	1.0000	\\
{\tt case24\_ieee\_rts}	&	63 333.39	&	63 335.67	&	1.0000	&	1.0000	&	0.9999	&	1.0000	&	1.0000	&	1.0000	&	1.0000	&	1.0000	\\
{\tt case30}	&	576.89	&	576.89	&	1.0000	&	1.0000	&	0.9993	&	0.9993	&	1.0000	&	1.0000	&	1.0000	&	1.0000	\\
{\tt case\_ieee30}	&	8 902.67	&	8 902.75	&	1.0000	&	1.0000	&	1.0000	&	1.0000	&	1.0000	&	1.0000	&	1.0000	&	1.0000	\\
{\tt case39}	&	41 850.24	&	41 852.45	&	1.0000	&	1.0000	&	1.0000	&	1.0000	&	1.0000	&	1.0000	&	1.0000	&	1.0000	\\
{\tt case57}	&	41 682.14	&	41 688.61	&	0.9997	&	1.0000	&	0.9996	&	0.9999	&	1.0000	&	1.0000	&	1.0000	&	1.0000	\\
{\tt case89pegase}	&	5 803.58	&	5 804.23	&	0.9996	&	1.0000	&	0.9992	&	0.9999	&	1.0001	&	1.0000	&	1.0001	&	1.0001	\\
\multicolumn{3}{l|}{\bf Average}	&	\bf 0.9998	&	\bf 1.0000	&	\bf 0.9996	&	\bf 0.9999	&	\bf 1.0000	&	\bf 1.0000	&	\bf 1.0000	&	\bf 1.0000	\\
\midrule
\multicolumn{11}{c}{\it Medium-scale instances}\\
\midrule
{\tt case118}	&	129 526.00	&	129 662.15	&	0.9987	&	1.0000	&	0.9987	&	0.9997	&	1.0003	&	1.0000	&	1.0003	&	1.0000	\\
{\tt case\_ACTIV\_SG\_200}	&	27 552.84	&	27 553.52	&	1.0000	&	1.0000	&	1.0000	&	1.0000	&	*1.0000	&	1.0000	&	*1.0000	&	1.0000	\\
{\tt case\_illinois200}	&	36 738.07	&	36 743.37	&	1.0000	&	1.0000	&	0.9998	&	0.9998	&	*1.0000	&	1.0000	&	1.0000	&	1.0000	\\
{\tt case300}	&	718 938.90	&	719 154.24	&	0.9994	&	1.0000	&	0.9989	&	0.9992	&	1.0002	&	1.0000	&	1.0002	&	1.0002	\\
{\tt case\_ACTIV\_SG\_500}	&	70 316.47	&	72 454.07	&	0.9875	&	1.0000	&	0.9751	&	0.9792	&	1.0016	&	1.0033	&	1.0000	&	1.0018	\\
\multicolumn{3}{l|}{\bf Average}	&	\bf 0.9950	&	\bf 1.0000	&	\bf 0.9902	&	\bf 0.9919	&	\bf 1.0007	&	\bf 1.0013	&	\bf 1.0001	&	\bf 1.0007	\\
\midrule
\multicolumn{11}{c}{\it Large-scale instances}\\
\midrule
{\tt case1354pegase}	&	73 999.72	&	74 005.77	&	0.9998	&	1.0000	&	0.9995	&	0.9998	&	1.0000	&	1.0000	&	1.0000	&	1.0000	\\
{\tt case2383wp}	&	1 856 849.35	&	1 862 626.70	&	0.9996	&	1.0000	&	0.9981	&	0.9992	&	1.0039	&	1.0000	&	1.0068	&	1.0033	\\
{\tt case2736sp}	&	1 306 771.84	&	1 307 134.47	&	0.9997	&	1.0000	&	0.9988	&	0.9998	&	1.0001	&	1.0000	&	1.0002	&	1.0001	\\
{\tt case2737sop}	&	777 020.01	&	777 337.82	&	0.9997	&	1.0000	&	0.9990	&	0.9997	&	1.0001	&	1.0000	&	1.0001	&	1.0001	\\
{\tt case2746wop}	&	1 206 874.60	&	1 207 872.71	&	0.9994	&	1.0000	&	0.9983	&	0.9994	&	1.0001	&	1.0000	&	1.0002	&	1.0001	\\
{\tt case2746wp}	&	1 630 617.28	&	1 631 018.94	&	0.9997	&	1.0000	&	0.9988	&	0.9996	&	1.0001	&	1.0000	&	1.0002	&	1.0000	\\
{\tt case2869pegase}	&	133 867.21	&	133 877.91	&	0.9998	&	1.0000	&	0.9994	&	0.9997	&	1.0000	&	1.0000	&	1.0000	&	1.0000	\\
{\tt case3012wp}	&	2 578 533.79	&	2 590 254.80	&	0.9987	&	1.0000	&	0.9968	&	0.9985	&	--	&	1.0000	&	1.0009	&	1.0013	\\
{\tt case3120sp}	&	2 135 933.93	&	2 143 654.67	&	--	&	1.0000	&	0.9975	&	0.9988	&	--	&	1.0002	&	--	&	1.0000	\\
{\tt case3375wp}	&	7 397 491.46	&	*7 409 223.47	&	0.9998	&	1.0000	&	0.9988	&	0.9994	&	*1.0002	&	*1.0000	&	*1.0002	&	*1.0004	\\
\multicolumn{3}{l|}{\bf Average}	&	\bf --	&	\bf 1.0000	&	\bf 0.9985	&	\bf 0.9994	&	\bf --	&	\bf 1.0000	&	\bf --	&	\bf 1.0005	\\
\bottomrule
\end{tabular}
\end{table*}

\begin{table*}[!t]
\footnotesize
\centering
\caption{Loss minimization: Normalized optimal values}\label{table2:loss1}
\begin{tabular}{@{}lrr|rrrr|rrrr@{}}
\toprule
Test case & $\lb{\upsilon}$ [MW] & $\ub{\upsilon}$ [MW] & \multicolumn{4}{c|}{Normalized $\hat{\upsilon}_R$} & \multicolumn{4}{c}{Normalized $\ub{\upsilon}_{SP}$} \\
\cmidrule(r){4-7} \cmidrule(l){8-11} &&& SDR1 & SDR2 & TCR1 & TCR2 & SDR1 & SDR2 &  TCR1 & TCR2 \\
\midrule
\multicolumn{11}{c}{\it Small-scale instances}\\
\midrule
{\tt case14}	&	259.49	&	259.49	&	1.0000	&	1.0000	&	1.0000	&	1.0000	&	1.0000	&	1.0000	&	1.0000	&	1.0000	\\
{\tt case24\_ieee\_rts}	&	2 875.33	&	2 875.37	&	1.0000	&	1.0000	&	0.9999	&	1.0000	&	1.0000	&	1.0000	&	1.0000	&	1.0000	\\
{\tt case30}	&	191.09	&	191.09	&	1.0000	&	1.0000	&	0.9999	&	0.9999	&	1.0000	&	1.0000	&	1.0000	&	1.0000	\\
{\tt case\_ieee30}	&	284.68	&	284.71	&	0.9997	&	1.0000	&	0.9997	&	1.0000	&	1.0000	&	1.0000	&	1.0000	&	1.0000	\\
{\tt case39}	&	6 283.20	&	6 283.44	&	1.0000	&	1.0000	&	1.0000	&	1.0000	&	1.0000	&	1.0000	&	1.0000	&	1.0000	\\
{\tt case57}	&	1 260.83	&	1 260.98	&	0.9998	&	1.0000	&	0.9998	&	1.0000	&	1.0000	&	1.0000	&	1.0000	&	1.0000	\\
{\tt case89pegase}	&	5 803.58	&	5 804.23	&	0.9996	&	1.0000	&	0.9992	&	0.9999	&	1.0001	&	1.0000	&	1.0001	&	1.0001	\\
\multicolumn{3}{l|}{\bf Average}	&	\bf 0.9998	&	\bf 1.0000	&	\bf 0.9997	&	\bf 1.0000	&	\bf 1.0000	&	\bf 1.0000	&	\bf 1.0000	&	\bf 1.0000	\\
\midrule
\multicolumn{11}{c}{\it Medium-scale instances}\\
\midrule
{\tt case118}	&	4 250.75	&	4 251.17	&	1.0000	&	1.0000	&	1.0000	&	1.0000	&	1.0000	&	1.0000	&	1.0000	&	1.0000	\\
{\tt case\_ACTIV\_SG\_200}	&	1 483.23	&	1 483.32	&	1.0000	&	1.0000	&	0.9999	&	0.9999	&	1.0000	&	1.0000	&	*1.0000	&	1.0000	\\
{\tt case\_illinois200}	&	2 245.96	&	2 246.22	&	1.0000	&	1.0000	&	0.9999	&	0.9999	&	*1.0000	&	1.0000	&	1.0000	&	1.0000	\\
{\tt case300}	&	23 723.13	&	23 725.42	&	0.9999	&	1.0000	&	0.9996	&	0.9997	&	1.0000	&	1.0001	&	1.0001	&	1.0001	\\
{\tt case\_ACTIV\_SG\_500}	&	7 815.64	&	7 817.25	&	1.0000	&	1.0000	&	0.9998	&	0.9998	&	1.0000	&	1.0000	&	1.0000	&	1.0000	\\
\multicolumn{3}{l|}{\bf Average}	&	\bf 1.0000	&	\bf 1.0000	&	\bf 0.9998	&	\bf 0.9999	&	\bf 1.0000	&	\bf 1.0000	&	\bf 1.0000	&	\bf 1.0000	\\
\midrule
\multicolumn{11}{c}{\it Large-scale instances}\\
\midrule
{\tt case1354pegase}	&	73 999.72	&	74 005.77	&	0.9998	&	1.0000	&	0.9995	&	0.9998	&	1.0000	&	1.0000	&	1.0000	&	1.0000	\\
{\tt case2383wp}	&	24 967.98	&	24 983.29	&	1.0000	&	1.0000	&	0.9992	&	0.9997	&	1.0000	&	1.0000	&	1.0001	&	1.0001	\\
{\tt case2736sp}	&	18 324.27	&	18 329.28	&	0.9999	&	1.0000	&	0.9993	&	0.9998	&	1.0001	&	1.0000	&	1.0001	&	1.0000	\\
{\tt case2737sop}	&	11 391.09	&	11 394.38	&	0.9997	&	1.0000	&	0.9992	&	0.9997	&	1.0001	&	1.0000	&	1.0001	&	1.0001	\\
{\tt case2746wop}	&	19 197.71	&	19 209.58	&	0.9996	&	1.0000	&	0.9989	&	0.9996	&	1.0001	&	1.0000	&	1.0001	&	1.0001	\\
{\tt case2746wp}	&	25 255.49	&	25 260.72	&	0.9999	&	1.0000	&	0.9992	&	0.9997	&	1.0001	&	1.0000	&	1.0001	&	1.0000	\\
{\tt case2869pegase}	&	133 867.21	&	133 877.91	&	0.9998	&	1.0000	&	0.9994	&	0.9997	&	1.0000	&	1.0000	&	1.0000	&	1.0000	\\
{\tt case3012wp}	&	27 611.41	&	27 648.22	&	0.9995	&	1.0000	&	0.9986	&	0.9994	&	1.0001	&	1.0000	&	1.0001	&	1.0003	\\
{\tt case3120sp}	&	21 463.22	&	21 524.52	&	0.9998	&	1.0000	&	0.9991	&	0.9997	&	--	&	1.0000	&	--	&	1.0000	\\
{\tt case3375wp}	&	48 950.55	&	49 002.59	&	0.9997	&	1.0000	&	0.9993	&	0.9995	&	--	&	1.0027	&	--	&	1.0000	\\
\multicolumn{3}{l|}{\bf Average}	&	\bf 0.9998	&	\bf 1.0000	&	\bf 0.9992	&	\bf 0.9996	&	\bf --	&	\bf 1.0003	&	\bf --	&	\bf 1.0001	\\
\bottomrule
\end{tabular}
\end{table*}

\begin{table}[!t]
\footnotesize
\centering
\caption{Cost minimization: Optimality gaps and computation times}\label{table2:cost2}
\resizebox{\linewidth}{!}{
\begin{tabular}{@{}l|rrrr|rrrr@{}}
\toprule
Test case & \multicolumn{4}{c|}{Optimality gap [\%]} & \multicolumn{4}{c}{Computation time [s]} \\
\cmidrule(r){2-5} \cmidrule(l){6-9} & SDR1 & SDR2 & TCR1 & TCR2 & SDR1 & SDR2 &  TCR1 & TCR2 \\
\midrule
\multicolumn{9}{c}{\it Small-scale instances}\\
\midrule
{\tt case14}	&	0.00	&	0.00	&	0.00	&	0.00	&	0.15	&	0.13	&	0.12	&	0.36	\\
{\tt case24\_ieee\_rts}	&	0.01	&	0.00	&	0.01	&	0.00	&	0.19	&	0.55	&	0.54	&	0.28	\\
{\tt case30}	&	0.00	&	0.00	&	0.07	&	0.07	&	0.12	&	0.12	&	0.20	&	0.20	\\
{\tt case\_ieee30}	&	0.00	&	0.00	&	0.00	&	0.00	&	0.12	&	0.11	&	0.13	&	0.14	\\
{\tt case39}	&	0.01	&	0.01	&	0.01	&	0.01	&	0.43	&	0.41	&	0.61	&	0.52	\\
{\tt case57}	&	0.05	&	0.02	&	0.05	&	0.03	&	0.47	&	0.21	&	0.20	&	0.24	\\
{\tt case89pegase}	&	0.06	&	0.01	&	0.10	&	0.03	&	2.26	&	1.64	&	0.76	&	0.72	\\
\bf Average	&	\bf 0.03	&	\bf 0.01	&	\bf 0.05	&	\bf 0.02	&	\bf 0.91	&	\bf 0.69	&	\bf 0.45	&	\bf 0.42	\\
\midrule
\multicolumn{9}{c}{\it Medium-scale instances}\\
\midrule
{\tt case118}	&	0.27	&	0.11	&	0.27	&	0.13	&	0.30	&	0.36	&	0.38	&	0.40	\\
{\tt case\_ACTIV\_SG\_200}	&	*0.00	&	0.00	&	*0.01	&	0.01	&	0.82	&	0.79	&	0.57	&	0.55	\\
{\tt case\_illinois200}	&	*0.02	&	0.01	&	0.03	&	0.03	&	0.94	&	1.41	&	0.73	&	0.75	\\
{\tt case300}	&	0.11	&	0.03	&	0.16	&	0.13	&	1.23	&	2.25	&	0.99	&	0.89	\\
{\tt case\_ACTIV\_SG\_500}	&	4.31	&	3.27	&	5.36	&	5.14	&	4.37	&	4.53	&	3.85	&	3.14	\\
\bf Average	&	\bf 1.69	&	\bf 1.26	&	\bf 2.10	&	\bf 2.00	&	\bf 2.23	&	\bf 2.60	&	\bf 1.92	&	\bf 1.63	\\
\midrule
\multicolumn{9}{c}{\it Large-scale instances}\\
\midrule
{\tt case1354pegase}	&	0.03	&	0.01	&	0.06	&	0.03	&	13.35	&	11.95	&	7.53	&	6.66	\\
{\tt case2383wp}	&	0.74	&	0.31	&	1.17	&	0.71	&	307.06	&	325.67	&	14.19	&	12.97	\\
{\tt case2736sp}	&	0.07	&	0.03	&	0.16	&	0.07	&	426.85	&	398.85	&	13.84	&	12.05	\\
{\tt case2737sop}	&	0.08	&	0.04	&	0.16	&	0.08	&	381.34	&	372.61	&	11.23	&	12.75	\\
{\tt case2746wop}	&	0.15	&	0.08	&	0.27	&	0.15	&	477.80	&	425.23	&	10.96	&	11.68	\\
{\tt case2746wp}	&	0.06	&	0.02	&	0.16	&	0.07	&	418.71	&	414.08	&	11.27	&	11.19	\\
{\tt case2869pegase}	&	0.03	&	0.01	&	0.07	&	0.04	&	55.71	&	80.07	&	21.27	&	19.78	\\
{\tt case3012wp}	&	--	&	0.45	&	0.87	&	0.73	&	1 670.60	&	1 618.99	&	12.27	&	12.17	\\
{\tt case3120sp}	&	--	&	0.38	&	--	&	0.48	&	763.19	&	672.87	&	13.64	&	13.22	\\
{\tt case3375wp}	&	*0.20	&	*0.16	&	*0.29	&	*0.24	&	1 638.91	&	1 763.90	&	20.47	&	18.89	\\
\bf Average	&	\bf --	&	\bf 0.16	&	\bf --	&	\bf 0.27	&	\bf 684.22	&	\bf 678.29	&	\bf 14.15	&	\bf 13.63	\\
\bottomrule
\end{tabular}
}
\end{table}

\begin{table}[!t]
\footnotesize
\centering
\caption{Loss minimization: Optimality gaps and computation times}\label{table2:loss2}
\resizebox{\linewidth}{!}{
\begin{tabular}{@{}l|rrrr|rrrr@{}}
\toprule
Test case & \multicolumn{4}{c|}{Optimality gap [\%]} & \multicolumn{4}{c}{Computation time [s]} \\
\cmidrule(r){2-5} \cmidrule(l){6-9} & SDR1 & SDR2 & TCR1 & TCR2 & SDR1 & SDR2 &  TCR1 & TCR2 \\
\midrule
\multicolumn{9}{c}{\it Small-scale instances}\\
\midrule
{\tt case14}	&	0.00	&	0.00	&	0.00	&	0.00	&	0.12	&	0.12	&	0.12	&	0.12	\\
{\tt case24\_ieee\_rts}	&	0.00	&	0.00	&	0.00	&	0.00	&	0.16	&	0.18	&	0.39	&	0.20	\\
{\tt case30}	&	0.00	&	0.00	&	0.01	&	0.01	&	0.13	&	0.13	&	0.17	&	0.17	\\
{\tt case\_ieee30}	&	0.04	&	0.01	&	0.05	&	0.02	&	0.16	&	0.34	&	0.49	&	0.28	\\
{\tt case39}	&	0.00	&	0.00	&	0.00	&	0.00	&	0.17	&	0.19	&	0.22	&	0.22	\\
{\tt case57}	&	0.03	&	0.01	&	0.03	&	0.01	&	0.19	&	0.23	&	0.18	&	0.49	\\
{\tt case89pegase}	&	0.06	&	0.01	&	0.10	&	0.03	&	2.26	&	1.64	&	0.76	&	0.72	\\
\bf Average	&	\bf 0.03	&	\bf 0.01	&	\bf 0.04	&	\bf 0.02	&	\bf 0.82	&	\bf 0.66	&	\bf 0.41	&	\bf 0.43	\\
\midrule
\multicolumn{9}{c}{\it Medium-scale instances}\\
\midrule
{\tt case118}	&	0.01	&	0.01	&	0.01	&	0.01	&	0.29	&	0.27	&	0.30	&	0.30	\\
{\tt case\_ACTIV\_SG\_200}	&	0.01	&	0.01	&	*0.01	&	0.01	&	0.73	&	0.74	&	0.48	&	1.06	\\
{\tt case\_illinois200}	&	*0.01	&	0.01	&	0.02	&	0.02	&	1.83	&	0.92	&	1.70	&	0.74	\\
{\tt case300}	&	0.02	&	0.02	&	0.06	&	0.05	&	1.02	&	0.91	&	0.98	&	0.88	\\
{\tt case\_ACTIV\_SG\_500}	&	0.03	&	0.02	&	0.04	&	0.04	&	5.82	&	4.56	&	4.59	&	3.95	\\
\bf Average	&	\bf 0.02	&	\bf 0.02	&	\bf 0.04	&	\bf 0.03	&	\bf 2.85	&	\bf 2.21	&	\bf 2.32	&	\bf 2.00	\\
\midrule
\multicolumn{9}{c}{\it Large-scale instances}\\
\midrule
{\tt case1354pegase}	&	0.03	&	0.01	&	0.06	&	0.03	&	13.35	&	11.95	&	7.53	&	6.66	\\
{\tt case2383wp}	&	0.07	&	0.06	&	0.15	&	0.10	&	252.55	&	271.06	&	13.94	&	12.56	\\
{\tt case2736sp}	&	0.05	&	0.03	&	0.11	&	0.05	&	434.61	&	280.45	&	10.52	&	9.58	\\
{\tt case2737sop}	&	0.06	&	0.03	&	0.12	&	0.06	&	210.99	&	215.05	&	8.67	&	8.55	\\
{\tt case2746wop}	&	0.10	&	0.06	&	0.18	&	0.12	&	394.81	&	444.52	&	9.88	&	11.20	\\
{\tt case2746wp}	&	0.04	&	0.02	&	0.11	&	0.05	&	431.06	&	403.43	&	9.70	&	9.72	\\
{\tt case2869pegase}	&	0.03	&	0.01	&	0.07	&	0.04	&	55.71	&	80.07	&	21.27	&	19.78	\\
{\tt case3012wp}	&	0.19	&	0.13	&	0.28	&	0.22	&	2 354.15	&	2 426.65	&	16.82	&	17.45	\\
{\tt case3120sp}	&	--	&	0.28	&	--	&	0.32	&	1 193.18	&	973.71	&	16.66	&	19.16	\\
{\tt case3375wp}	&	--	&	0.38	&	--	&	0.15	&	1 415.74	&	1 588.19	&	28.25	&	39.33	\\
\bf Average	&	\bf --	&	\bf 0.12	&	\bf --	&	\bf 0.12	&	\bf 753.59	&	\bf 749.07	&	\bf 15.09	&	\bf 16.54	\\
\bottomrule
\end{tabular}
}
\end{table}


%

%


\ifCLASSOPTIONcaptionsoff
  \newpage
\fi



%
\bibliographystyle{unsrt}
\bibliography{../drmth_biblio}
%
%

%

\begin{IEEEbiographynophoto}{Christian Bingane}
(S'18) received the B.Eng. degree in electrical engineering in 2014 from Polytechnique Montreal, Montreal, QC, Canada, where he is currently working toward the Ph.D. degree in applied mathematics. He is currently a student member of the GERAD research center.

His research interests include optimization in power systems and conic programming. He is concerned with using linear programming, second-order cone programming or semidefinite programming to provide guaranteed global optimal solution to the optimal power flow problem for a large-scale power system.
\end{IEEEbiographynophoto}

\begin{IEEEbiographynophoto}{Miguel F. Anjos}
(M'07--SM'18) received the B.Sc., the M.S. and the Ph.D. degrees from McGill University, Montreal, QC, Canada; Stanford University, Stanford, CA, USA; and the University of Waterloo, Waterloo, ON, Canada in 1992, 1994 and 2001 respectively. 

He is currently Chair of Operational Research at the School of Mathematics, University of Edinburgh, Edinburgh, Scotland, UK, and Professor with the Department of Mathematics and Industrial Engineering, Polytechnique Montreal, Montreal, QC, Canada, where he holds the NSERC-Hydro-Quebec-Schneider Electric Industrial Research Chair, and an Inria International Chair. He is a Licensed Professional Engineer in Ontario, Canada. He served for five years as Editor-in-Chief of Optimization and Engineering, and serves on several editorial boards.

His allocades include a Canada Research Chair, the M\'eritas Teaching Award, a Humboldt Research Fellowship, the title of EUROPT Fellow, and the Queen Elizabeth II Diamond Jubilee Medal. He is an elected Fellow of the Canadian Academy of Engineering.
\end{IEEEbiographynophoto}


\begin{IEEEbiographynophoto}{S\'ebastien Le~Digabel}
received the M.Sc.A. and the Ph.D. degrees in applied mathematics from Polytechnique Montreal, Montreal, Quebec, Canada in 2002 and 2008 respectively. He was a postdoctoral fellow with the IBM Watson Research Center and the University of Chicago in 2010 and 2011.

He is currently an Associate Professor with the Department of Mathematics and Industrial Engineering, Polytechnique Montreal, Montreal, QC, Canada, and a regular member of the GERAD research center.

His research interests include the analysis and development of algorithms for derivative-free and blackbox optimization, and the design of related software. All of his work on derivative-free optimization is included in the NOMAD software, a free package for blackbox optimization available at {\tt www.gerad.ca/nomad}.
\end{IEEEbiographynophoto}




\end{document}